\documentclass[12pt,twoside,a4paper]{amsart}
\usepackage{amssymb}
\date{\today}

\def\deg{\text{deg}\,}
\def\sr{\stackrel}
\def\lsigma{\text{\Large $\sigma$}}

\def\w{\wedge}

\def\dbar{\bar\partial}

\def\C{{\mathbb C}}
\def\w{{\wedge}}

\def\B{{\mathbb B}}

\def\Hom{{\rm Hom\, }}
\def\codim{{\rm codim\,}}

\def\Re{{\rm Re\,  }}

\def\be{\begin{equation}}
\def\ee{\end{equation}}

\newtheorem{thm}{Theorem}[section]
\newtheorem{lma}[thm]{Lemma}
\newtheorem{cor}[thm]{Corollary}
\newtheorem{prop}[thm]{Proposition}

\theoremstyle{definition}

\theoremstyle{remark}

\newtheorem{preremark}{Remark}
\newtheorem{preex}{Example}

\newenvironment{remark}{\begin{preremark}}{\qed\end{preremark}}
\newenvironment{ex}{\begin{preex}}{\qed\end{preex}}

\numberwithin{equation}{section}

\title[Explicit versions of the Brian\c con-Skoda theorem \ldots]
{Explicit versions of the Brian\c con-Skoda theorem with variations}

\begin{document}

\date{\today}

\author{Mats Andersson}

\address{Department of Mathematics\\Chalmers University of Technology and the University of G\"oteborg\\S-412 96 G\"OTEBORG\\SWEDEN}

\email{matsa@math.chalmers.se}



\thanks{The author was
  partially supported by the Swedish 
  Research Council}


\maketitle

\section{Introduction}

Let $\phi,f_1,\ldots,f_m$ be holomorphic functions in a neighborhood of the origin 
in $\C^n$. The Brian\c con-Skoda theorem, \cite{BS}, states that
 $\phi^{\min(n,m)}$ belongs to the ideal $(f)$ generated by $f_j$ if
 $|\phi|\le C|f|$.
This  condition is equivalent to that $\phi$ belongs to the integral closure
of the ideal $(f)$. The original proof is based on  Skoda's $L^2$-estimates
in \cite{Sk1}, see Remark~\ref{rustan} below, 
and actually gives the stronger statement that 
$\phi\in(f)$ if $|\phi|\le C|f|^{\min(n,m)}$. An explicit proof based on Berndtsson's
division formula,  \cite{BB},   and multivariable residue calculus appeared in
\cite{BGVY}, see also \cite{Elk} for a special case. 
There are purely algebraic versions in 
more arbitrary rings due to  Lipman and Teissier, \cite{LT}.

In general this result cannot be improved but for certain tuples $f_j$
a much weaker size condition on $\phi$ is enough to guarantee that
$\phi$ belongs to $(f)$. For instance, the ideal $(f)^2$ is generated
by the $m(m+1)/2$ functions  $g_{jk}=f_jf_k$, and $|f|^2\sim |g|$, so if we apply the
previous result we get that 
$\phi\in(f)^2$ if  $|\phi|\le C|f|^{\min(2n,m(m+1))}$. However, in this case
actually the power $\min(n,m)+1$ is enough. In general we have

\begin{thm}[Brian\c con-Skoda]\label{bs}
If $f=(f_1,\ldots, f_m)$ and $\phi$ are holomorphic  at $0$ in $\C^n$ and 
$|\phi|\le C|f|^{\min(m,n)+r-1}$, then $\phi\in(f)^r$.
\end{thm}

This more general formulation follows in a similar way as the case $r=1$ by 
$L^2$-methods as well as by a (small modification of)  the argument in \cite{BGVY}.
In \cite{A2} we gave a   somewhat different proof of  the case $r=1$
by means of residue calculus and in this note we extend this method
to get various related for product ideals, as well as
 the general case of  Theorem~\ref{bs}.
In the first one we consider several possibly  different tuples.

\begin{thm}\label{ma1}
Let $f_j$, $j=1,\ldots,r$,  be $m_j$-tuples of holomorphic functions at $0\in\C^n$
and assume that 
$$
|\phi|\le C|f_1|^{s_1}\cdots |f_r|^{s_r}
$$
for all $s$ such that $s_1+\cdots +s_r\le n+r-1$ and $1\le s_j\le m_j$.
Then $\phi\in (f_1)\cdots (f_r)$.
\end{thm}

Notice that this immediately implies Theorem~\ref{bs} in the case $m\ge n$ by just choosing
all $f_j=f$.
In certain cases
Theorem~\ref{ma1} can be improved, as one can  see by taking $f_j=f$ and $m<n$
and compare with   Theorem~\ref{bs}. Another case is when all the functions
in the various tuples $f_j$ together form a regular sequence.

\begin{thm}\label{ma2}
Let $f_j$, $j=1,\ldots,m$,  be $m_j$-tuples of holomorphic functions at $0\in\C^n$
and assume that the codimension of $\{f_1=\cdots =f_r=0\}$ is
$m_1+\cdots +m_r$. If 
$$
|\phi|\le C \min(|f_1|^{m_1},\ldots,|f_r|^{m_r}),
$$
then $\phi\in (f_1)\cdots (f_r)$.
\end{thm}

We have not seen these two latter results in the literature although they might belong to
 the folklore. 
In the algebraic setting there are several results related to the 
Brian\c con-Skoda (and Lipman-Teissier) theorem, see, e.g.,
\cite{Te} and \cite{HH} and the references given there.

\begin{remark}\label{rustan} 
The  Brian\c con-Skoda theorem  follows  by direct  a  application  of Skoda's
$L^2$-estimate, \cite{Sk1} and \cite{Sk2},  if $m\le n$. In fact,  if $\psi$ is any  plurisubharmonic function,
the $L^2$-estimate guarantees a holomorphic solution to
$f\cdot u=\phi$ such that
$$
\int_{X\setminus Z}\frac{|u|^2}{|f|^{2(\min(m,n+1)-1+\epsilon)}}e^{-\psi}dV<\infty
$$
provided that
$$
\int_{X\setminus Z}\frac{|\phi|^2}{|f|^{2(\min(m,n+1)+\epsilon)}}e^{-\psi}dV<\infty.
$$
If $|\phi|\le C|f|^{m}$, the second integral is finite (taking  $\psi=0$)
if $\epsilon$ is small enough, and thus Skoda's theorem provides the desired solution.
The case when $r>1$ is  obtained by iteration. 
If $m>n$   a direct use of  the $L^2$-estimate will not give the
desired result. However, see \cite{Dem}, in this case 
one can find an  $n$-tuple $\tilde f$ such that $(\tilde f)\subset (f)$ and 
$|\tilde f|\sim |f|$, and the theorem then follows by applying the $L^2$-estimate
to the tuple $f'$.

In the same way, Theorem~\ref{ma1} can easily be proved  from the $L^2$-
estimate if $m_1+\cdots +m_r\le n+r-1$. 
To see this, assume for simplicity that $r=2$,
and that $|\phi|\le C|f_1|^{m_1}|f_2|^{m_2}$.  Choosing
$\psi=2(m_1+\epsilon)\log|f_1|$, Skoda's theorem gives 
a solution to  $f_2\cdot u=\phi$ such that 
$$
\int_{X\setminus Z}\frac{|u|^2}{|f_1|^{2(m_1+\epsilon)}}
dV<\infty.
$$
Another application then gives $v_j$ such that $f_1 \cdot v_j=u_j$. This means that
$\phi$ belongs to $(f_1)(f_2)$.
However, we do not know   whether  one can derive Theorem~\ref{ma2}
from the $L^2$-estimate when $m_1+\cdots +m_r>n+r-1$.
\end{remark}

Now consider an  $r\times m$ matrix $f_j^k$ of holomorphic functions, $r\le m$, with rows
$f_1,\ldots, f_r$. We let $F$ be the $m!/(m-r)!r!$ tuple of functions
$\det(f_j^{I_k})$ for increasing multiindices $I$ of length $r$.
We will refer to $F$ as the determinant 
of $f$. If $f_j$ are the rows of the matrix, considered as sections of the trivial bundle
$E^*$, then $F$ is just the section $f_r\w\ldots \w f_1$ of the bundle $\Lambda^r E^*$.
Our next result is a Brian\c con-Skoda type result for
the tuple $F$. It turns out that
it is enough with a much less power than $m!/(m-r)!r!$.
Let $Z$ be the zero set of $F$ and notice that  $\codim Z\le m-r+1$;
this is  easily seen by Gauss elimination. 

\begin{thm}\label{ma3}
Let $F$ be the determinant of the holomorphic matrix $f$ as above. If 
$$
|\phi|\le C|F|^{\min(n,m-r+1)},
$$
then $\phi\in(F)$.
\end{thm}

\begin{remark}\label{sporr}
This result is closely related to the following statement which was proved in
\cite{A4}.
{\it Suppose that $\phi$ is an $r$-tuple of holomorphic functions and let
$\|\phi\|$ be the pointwise norm induced by $f$, i.e., 
$\|\phi\|=\det(ff^*)\langle (ff^*)^{-1}\phi,\phi\rangle$. If
$$
\|\phi\|\lesssim |F|^{\min(n,m-r+1)},
$$
then $f\psi=\phi$ has a local holomorphic solution. }
\end{remark}

\begin{remark}
Another related situation is when $f$ is a section of a bundle $E^*$, 
$\phi$ takes values in $\Lambda^\ell E$,  and we ask for a holomorphic section
$\psi$ of $\Lambda^{\ell+1}E$ such that $\delta_f\psi=\phi$, provided that
the necessary compatibility condition $\delta_f\phi=0$ is fulfilled.
Let $p=\codim\{f=0\}$. Then a sufficient condition  is that
$$
|\phi|\le  C|f|^{\min(n,m-\ell)}
$$
if $\ell\le m-p$, whereas there is no condition at all if $\ell> m-p$,
see Theorems~1.2 and 1.4 and Corollary~1.5 in \cite{A2}. 
\end{remark}

Theorem~\ref{ma3}  is proved by constructing a certain residue current $R$ with support on the
analytic set $Z$,  such that
$R\phi=0$ implies that $\phi$ belongs to the ideal $(F)$ locally.   The size conditions
of $\phi$ then implies that $R\phi=0$ by brutal force,
see Theorem~\ref{hase} below.   There may be 
more subtle reasons for annihilation.  For instance, in the generic case, i.e., when $\codim Z=m-r+1$,
even the converse statement holds; if $\phi$ is in the ideal $(F)$ 
then actually $R\phi=0$, see Theorem~\ref{hase}~(iv).
The analogous statement also holds  for the equation $f\psi=\phi$ 
in Remark~\ref{sporr},  see \cite{A4}.
These results are therefore  extensions of the
well-known duality theorem of Dickenstein-Sessa and Passare, \cite{DS} and
\cite{P1},  stating  that if  $f$ is a tuple that defines  a complete intersection, i.e.,
$\codim \{f=0\}=m$, then $\phi\in(f)$  if and only if
$\phi$ annihilates the Coleff-Herrera current defined by $f$.
For the analysis of the  residue current $R$ we use the  basic tools 
developed in \cite{BY}, \cite{BGVY}, \cite{PTY}, and \cite{BY3},   i.e.,
resolution of singularites by Hironaka's theorem followed by a toric resolution. 
Theorems~\ref{ma1} and \ref{ma2} (as well as Theorem~\ref{bs})
are obtained along the same lines, in  Section~\ref{poo}, by analysis of special choices of
matrix $f$. It might happen that
there are some similarities  with the methods used  here and
the algebraic methods introduced in \cite{EL}. 

\smallskip

By means of  the new  construction in \cite{A1} of division formulas
we get,    for a given holomorphic function $\phi$, a holomorphic decomposition
\begin{equation}\label{deco}
\phi=T\phi +S\phi,
\end{equation}
 such that   $T\phi$ belongs to the determinant ideal $(F)$,
and  $S\phi$ vanishes as soon as $\phi$ annihilates the residue current $R$.
In particular, this gives an explicit proof of Theorem~\ref{ma3}, and it  
leads to explicit proofs of Theorems~\ref{bs} to \ref{ma2} as well.

\smallskip
{\bf Acknowledgement }  I am grateful to Alain Yger as well as the referee
for several important remarks on a previous version.

\section{The ideal generated by the determinant section}\label{curr}
Although we are mainly interested in local results in this paper it is
 convenient to adopt  an invariant perspective. We therefore assume that
we have Hermitian vector bundles $E$ and $Q$ of ranks $m$ and $r\le m$,
respectively, 
over a complex $n$-dimensional manifold $X$, and a holomorphic morphism
$f\colon E\to Q$. We also assume that $f$ is generically surjective, i.e.,
that the analytic set $Z$ where $f$  is not surjective has at least
codimension $1$. If $\epsilon_j$ is a local holomorphic frame for
$Q$, then $f=f_1\otimes\epsilon_1+\cdots +f_r\otimes\epsilon_r$,
where $f_j$ are sections of the dual bundle $E^*$. Moreover,
$F=f_r\w\ldots \w f_1\otimes \epsilon_1\w\ldots\w\epsilon_r$
is an invariantly defined section of $\Lambda^r E^*\otimes\det Q^*$
that we will call the determinant section associated with  $f$. 
Notice that if $e_j$ is a local frame for $E$ with dual  frame
$e_j^*$ for  $E^*$, then
$f_j=\sum_1^m f_j^k e_k^*$, and
$$
F=\sum'_{|I|=r} F_I e^*_{I_1}\w\ldots\w e^*_{I_r},
$$
where the sum runs over increasing multiindices $I$ and
$
F_I=\det(F_j^{I_k}).
$
Let $S^\ell Q^*$ be the subbundle of $(Q^*)^{\otimes \ell}$ consisting of
symmetric tensors.  We introduce  the complex
\begin{multline}\label{komplex}
\cdots \sr{\delta_f} {\rightarrow}\Lambda^{r+k-1}E\otimes S^{k-1}Q^*\otimes\det Q^*\sr{\delta_f} {\rightarrow}\cdots \\
\sr{\delta_f} {\rightarrow}  \Lambda^{r+1}E\otimes Q^*\otimes\det Q^*\sr{\delta_f} {\rightarrow}
\Lambda^r E\otimes\det Q^*\sr{\delta_F} {\rightarrow} \C\to 0,
\end{multline}
where
$$
\delta_f=\sum_j \delta_{f_j}\otimes\delta_{\epsilon_j},
$$
$\delta_{f_j}$ and $\delta_{\epsilon_j}$ denote   interior multiplication on
$\Lambda E$ and from the left on $SQ^*\otimes\det Q^*$,  respectively, and
$$
\delta_F=\delta_f^r/r!=\delta_{f_r}\cdots\delta_{f_1}\otimes\delta_{\epsilon_1}\cdots\delta_{\epsilon_r}.
$$
It is readily checked that \eqref{komplex} actually is a complex. 
Notice that if $r=1$,  then \eqref{komplex}  is the usual Koszul complex.

\smallskip
In $X\setminus Z$ we let $\sigma_j$ be  the sections of $E$ with minimal norms such that
$f_k\sigma_j=\delta_{jk}$. Then $\sigma=\sigma_1\otimes\epsilon_1^*+\ldots +\sigma_r\otimes\epsilon_r^*$
is the section of $\Hom(Q,E)$ such that, for each section $\phi$ of $Q$,
 $v=\sigma\phi$ is the solution to $fv=\phi$ with pointwise minimal norm.  
We also have the  invariantly defined section
$$
\lsigma=\sigma_1\w\ldots\w\sigma_r\otimes\epsilon_r^*\w\ldots\w\epsilon_1^*
$$
of $\Lambda E\otimes =\det Q^*$, which in fact is  the section with minimal
norm such that $F\lsigma=1$, see, e.g., \cite{A4}.

\begin{ex}
Assume that  $E$ and $Q$ are trivial and let  $\epsilon_j$ be an ON-frame for $Q$ and
$e_j$ an ON-frame for $E$, with dual frame $e_j^*$. 
If $F=\sum_{|I|=r}' F_I e^*_{I_1}\w\ldots\w e_{I_r}^*$
as above, then
$$
\lsigma=\sum_{|I|=r}'\frac{\bar F_I}{|F|^2}e_{I_1}\w\ldots\w e_{I_r}.
$$
\end{ex}

We will   consider $(0,q)$-forms with values in  $\Lambda^{r+k-1}E\otimes S^{k-1} Q^*\otimes\det Q^*$, 
and it is convenient to consider them as sections of
$\Lambda^{r+k+q-1}(E\oplus T^*_{0,1}(X))\otimes S^{k-1}Q^*\otimes\det Q^*$, so that $\delta_f$ anti-commutes with
$\dbar$, and $\delta_F\dbar=(-1)^r \delta_F\dbar$. 
In what follows we let $\otimes$ denote usual tensor product all  $Q^*$-factors,
and wedge product of  $\Lambda (E\oplus T^*_{0,1}(X))$-factors.
Thus for instance 
$$
\sigma\otimes\lsigma=(\sum_1^r \sigma_j\otimes\epsilon^*_j)\otimes (\sigma_1\w\ldots\w\sigma_r\otimes
\epsilon_1^*\w\ldots\w\epsilon^*_r)=0.
$$
Moreover, for each $k\ge 1$, $(\dbar\sigma)^{\otimes(k-1)}$ is a symmetric tensor; more
precisely,
\begin{equation}\label{snorre}
(\dbar\sigma)^{\otimes(k-1)}=
\sum_{|\alpha|=k-1}(\dbar\sigma_1)^{\alpha_1}
\w\ldots\w (\dbar\sigma_r)^{\alpha_r}\otimes\epsilon_\alpha^*,
\end{equation}
where
$
\epsilon^*_\alpha=(\epsilon_1^*)^{\alpha_1}\dot\otimes\cdots\dot\otimes(\epsilon_r^*)^{\alpha_r}/
\alpha_1!\cdots \alpha_r!,
$
and   $\dot\otimes$ denotes symmetric tensor product. 
For each $k\ge 1$ we define in   $X\setminus Z$  the $(0,k-1)$-forms 
\begin{equation}\label{udef0}
u_k=(\dbar\sigma)^{\otimes(k-1)}\otimes\lsigma= 
\sigma_1\w\ldots\w \sigma_r\w(\dbar\sigma)^{\otimes(k-1)}\otimes \epsilon^*
\end{equation}
(where 
$\epsilon^*=\epsilon_r^*\w\ldots\w\epsilon_1^*$),
with values in   $\Lambda^{r+k-1}E\otimes S^{k-1}Q^*\otimes\det Q^*$.

\begin{prop}
In $X\setminus Z$ we have that 
\begin{equation} \label{gamma}
\delta_Fu_1=1, \quad \delta_f u_{k+1}=\dbar u_k, \ k\ge 1.
\end{equation}
\end{prop}

\begin{proof}
Since  $\delta_{\epsilon_j}$ act from the left,
and $\delta_{f_j}\dbar\sigma_\ell =0$ for all $\ell$, we have that
\begin{multline*}
\delta_f u_{k+1}=\delta_f\big[\sigma_1\w\ldots\w\sigma_r\w (\dbar\sigma)^{\otimes k}\otimes\epsilon^*\big] =\\
\delta_f\big[\sigma_1\w\ldots\w\sigma_r\w\dbar\sigma\big]\otimes(\dbar\sigma)^{\otimes(k-1)}\otimes\epsilon^* =\\
\sum_{j=1}^r\delta_{f_j}(\sigma_1\w\ldots\w\sigma_r)\w\dbar\sigma_j\otimes(\dbar\sigma)^{\otimes(k-1)}\otimes\epsilon^*=\\
\dbar(\sigma_1\w\ldots\w\sigma_r)\w(\dbar\sigma)^{\otimes(k-1)}\otimes\epsilon^*=
\dbar u_k.
\end{multline*}
Since $\delta_F u_1=F\lsigma=1$, the proposition is proved.
\end{proof}

If we let $u=u_1+u_2 +\cdots$, and let $\delta$ denote either $\delta_f$ or $\delta_F$,
then \eqref{gamma} can be written as  $(\delta-\dbar)u=1$.
To analyze the singularities of $u$ at $Z$ we will use the following lemma
(Lemma~4.1) from \cite{A4}.

\begin{lma}\label{fakt}
If   $F=F_0 F'$ for some holomorphic  function $F_0$ and 
non-vanishing holomorphic section $F'$, then 
$$
s'=F_0\sigma, \quad S'=F_0\lsigma
$$
are smooth across  $Z$.
\end{lma}

Notice that $|F|^{2\lambda}u$ and $\dbar|F|^{2\lambda}\w u$ are well-defined
forms in $X$ for $\Re\lambda>>0$.

\begin{thm}\label{hase}
\noindent (i)  The forms  $|F|^{2\lambda}u$ and $\dbar|F|^{2\lambda}\w u$ have analytic continuations
as currents in $X$ to $\Re\lambda>-\epsilon$. If
$U=|F|^{2\lambda}u|_{\lambda=0}$ and $R=\dbar|F|^{2\lambda}\w u|_{\lambda=0}$, then
$$
(\delta-\dbar)U=1-R.
$$

\smallskip

\noindent (ii)  The current $R$ has support on $Z$ and
$R=R_p+\cdots +R_\mu$,  where $p=\codim Z$ and $\mu=\min(n,m-r+1)$.
\smallskip

\noindent(iii) If $\phi$ is a holomorphic function and $R\phi=0$, then
locally $F\Psi=\phi$  has holomorphic solutions.
\smallskip

\noindent(iv)  If $\codim Z=m-r+1$ and $F\Psi =\phi$ has  a holomorphic solution, then
$R\phi=R_{m-r+1}\phi=0$.

\smallskip
\noindent(v)  If $|\phi|\le C |F|^\mu$, then $R\phi=0$.
\end{thm}

Here, of course,   $R_k=\dbar|F|^{2\lambda}\w u_k|_{\lambda=0}$ is the component of $R$
which is a $(0,k)$-current with values in $\Lambda^{r+k-1}E\otimes S^{k-1}Q^*\otimes\det Q^*$.

\begin{proof}
In the case $r=1$,  this  theorem is contained in Theorems~1.1 to 1.4 in \cite{A2},
and most parts of the proof are completely analogous. 
Therefore we just point out the necessary modifications.
By Hironaka's theorem and a further toric resolution,
following the technique developed in \cite{BY} and \cite{PTY},  we may assume that
locally $F=F_0F'$ as in  Lemma~\ref{fakt}. Since
moreover $\sigma\otimes\lsigma=0$, we have then that locally in the resolution 
$$
u_k=\frac{(\dbar s')^{\otimes(k-1)}\otimes S'}{F_0^k}.
$$
It is then easy to see that the proposed analytic extensions exist and we have that 
\begin{equation}\label{using}
U_k=\Big[\frac{1}{F_0^k}\Big](\dbar s')^{\otimes(k-1)}\otimes S',
\end{equation}
and
\begin{equation}\label{rgrad}
R_k=\dbar \Big[\frac{1}{F_0^k}\Big]\w(\dbar s')^{\otimes(k-1)}\otimes S',
\end{equation}
where $[1/F_0^k]$ is the usual principal value current. 
If $R\phi=0$, then $(\delta-\dbar)U\phi=\phi$, and hence by successively
solving the $\dbar$-equations
$
\dbar w_k=U_k\phi+\delta w_{k+1},
$
we finally get the holomorphic solution $\Psi=U_1\phi+\delta w_2$.
All parts but $(iv)$ now follow in a similar way as in \cite{A2}.
Notice in particular, that $k\le\min(n,m-r+1)$ in \eqref{rgrad}
for degree reasons, 
so that $R\phi=0$ if the hypothesis in (v) is satisfied.
As for $(iv)$,  let us  assume that we have a holomorphic  section $\Psi$ of
$\Lambda^rE\otimes\det Q^*$ such that $F\Psi=\phi$. 
If $\Psi=\psi\otimes \epsilon_*$, then
$F\Psi=\delta_{f_r}\cdots\delta_{f_1}\psi$.
Since $u_{m-r+1}$ has full degree in $e_j$ we have that
\begin{multline*}
u_{m-r+1}\phi= \phi \sigma_1\w\ldots\w\sigma_r\w(\dbar\sigma)^{\otimes(m-r)}
\otimes\epsilon^*=\\
(\delta_{f_r}\cdots\delta_{f_1}\psi)
\sigma_1\w\ldots\w\sigma_r\w(\dbar\sigma)^{\otimes(m-r)}
\otimes\epsilon^* =\\
\psi\w(\dbar\sigma)^{\otimes(m-r)}\otimes\epsilon^*|_{\lambda=0}=\\
(\dbar\sigma)^{\otimes(m-r)}\otimes\Psi=\dbar \big(\sigma\otimes(\dbar\sigma)^{\otimes(m-r+1)}\big)\otimes\Psi
=\dbar u'_{m-r}\otimes\Psi.
\end{multline*}
Since $\codim Z=m-r+1$ we have that $R=R_{m-r+1}$ according to part~(ii), so 
$$
R\phi= R_{m-r+1}\phi=\dbar|F|^{2\lambda}\w u_{m-r+1}\phi|_{\lambda=0}=
- \dbar \big(\dbar|F|^{2\lambda}\w u'_{m-r}\otimes\Psi|_{\lambda=0}\big).
$$
However, 
$$
\dbar|F|^{2\lambda}\w u'_{m-r}\otimes\Psi|_{\lambda=0}
$$
vanishes for degree reasons, precisely in the same way as $R_k$ vanishes for $k\le m-r$.
\end{proof}

\begin{proof}[Proof of  Theorem~\ref{ma3}]
If we consider the matrix $f$ as a morphism $E\to Q$,
for trivial bundles $E$ and $Q$,  the theorem  immediately follows from parts (v) and (iii)
of Theorem~\ref{hase}.
\end{proof}

\begin{remark}
As we have seen, the reason for the power $m-r+1$ in Theorem~\ref{ma3} (and in
part (v) of Theorem~\ref{hase}) when $n$ is large, is that the complex \eqref{komplex} terminates
at $k=m-r+1$. 
If one  tries to analyze the section $F$ by means of the usual Koszul complex  with respect to the
basis $(e_I)'_{|I|=r}$, then  one could hope  that for some miraculous reason
the corresponding forms $u_k$  would vanish when $k>m-r+1$, although
one has   $m!/(m-r)!r!$ dimensions (basis elements). However, this is not the case
in general. Take for instance the simplest  non-trivial case, $m=3$ and $r=2$, and
choose $f_1=(1,0,\xi_1)$, $f_2=(0,1,\xi_2)$ and choose the trivial metric.
Then $F_{12}=1,\  F_{13}=\xi_2,\  F_{23}=\xi_1,$ and 
$\lsigma={\bar F}/{|F|^2}$,   so that
$$
\lsigma_{12}=\frac{1}{1+|\xi_1|^2+|\xi_2|^2},
\lsigma_{13}=\frac{\bar\xi_2}{1+|\xi_1|^2+|\xi_2|^2},\quad \lsigma_{23}=\frac{\bar\xi_1}
{1+|\xi_1|^2+|\xi_2|^2}.
$$
Now  $m-r+1=2$, but if we form the usual Koszul complex, with say that basis
$\epsilon_1,\epsilon_2,\epsilon_3$, so that 
$$
\lsigma=\lsigma_{12}\epsilon_1+\lsigma_{13}\epsilon_2+\lsigma_{23}\epsilon_3=
\frac{1}{|F|^2}(\epsilon_1+\bar\xi_2\epsilon_2+\bar\xi_1\epsilon_3),
$$
we have 
$$
\lsigma\w(\dbar \lsigma)^2=\frac{2}{|F|^6}d\bar\xi_1\w d\bar\xi_2 \w \epsilon_1\w\epsilon_2\w\epsilon_3,
$$
and this  form is not zero. To get an example where
$Z$ is non-empty, one can   multiply $f$ with a   function $f_0$.
\end{remark}

\section{Products of ideals}\label{poo}

For $j=1,\dots,r$, let $E_j\to X$ be a Hermitian vector bundle of rank $m_j$ and
let $f_j$ be a section of $E_j^*$. Moreover, let $E=\oplus_1^r E_j$ and
let $Q\simeq \C^r$ with $ON$-basis $\epsilon_1,\ldots,\epsilon_r$. 
If we consider $f_j$ as sections of $E$, then $f=\sum_1^r f_j\otimes \epsilon_j$
is a morphism $E\to Q$. Moreover,   $F\Psi=\phi$ with $\Psi=\psi\otimes\epsilon^*$ as before, means that
$\delta_{f_r}\cdots\delta_{f_1}\psi=\phi$,  and hence that 
$\phi$ belongs to the product ideal $(f_1)\cdots (f_r)$.
To obtain such a solution $\Psi$ we proceed as in the previous section.  Notice that 
now $\sigma_j$ can be identified with the section  of  $E_j$ with minimal norm such that
$f_j\sigma_j=1$. Moreover, $|F|=|f_1|\cdots |f_r|$. 
In this case we therefore  have 
\begin{multline*}
R_k=\dbar|F|^{2\lambda}\w u_k=
\dbar(|f_1|^{2\lambda}\cdots |f_r|^{2\lambda})\w\sigma_1\w\ldots\w\sigma_r\w \\
\sum_{|\alpha|=k-1}(\dbar\sigma_1)^{\alpha_1}\w\ldots\w(\dbar\sigma_r)^{\alpha_r}
\otimes\epsilon_\alpha^*\otimes\epsilon^*|_{\lambda=0}.
\end{multline*}
For degree reasons $R_k$ will vanish unless
\begin{equation}\label{snark}
0\le \alpha_j\le m_j-1  \quad \text{and}\quad
\alpha_1+\cdots +\alpha_r\le n-1.
\end{equation}

\begin{proof}[Proof of Theorem~\ref{ma1}]
Consider the tuples $f_j$ as sections of $E_j$. 
For each $j$, let $e_{ji}$, $i=1,\ldots,m_j$,  be a local frame for $E_j$ so that
$f_j=\sum_{i=1}^{m_j} f_j^i e_{ji}^*$.
After a  suitable resolution as before we may assume that $f_1=f_1^0f_1'$, where
$f_1^0$ is holomorphic and $f_1'$ is a nonvanishing section of $E_1^*$.  After a further resolution
we may assume as well that $f_2=f_2^0f_2'$ etc.
Finally, therefore, we may assume that,  for each $j$, 
$f_j=f_j^0 f_j'$, where $f_j^0$ is holomorphic, and $f_j'$ is a non-vanishing section of
$E_j^*$. 
Therefore, $R_k$ is a sum of terms like
$$
\dbar(|f^0_1|^{2\lambda}\cdots |f^0_r|^{2\lambda}v^\lambda)\w 
\frac{\beta}{(f_1^0)^{\alpha_1+1}\cdots(f_r^0)^{\alpha_r+1}}\Big|_{\lambda=0},
$$
where $v$ is smooth and non-vanishing.
By the same argument as before this current is  annihilated by $\phi$ if
$|\phi|\le C|f_1|^{\alpha_1+1}\cdots |f_r|^{\alpha_r+1}$,
and in view of  \eqref{snark} and the hypothesis in the theorem, 
taking $s_j=\alpha_j+1$,
therefore $\phi$ annihilates  $R$. 
It now follows from Theorem~\ref{hase}~(iii) that
$F\Psi=\phi$ has a holomorphic solution, and thus  $\psi\in(f_1)\cdots(f_r)$.
\end{proof}

We can also easily obtain the Brian\c con-Skoda theorem.

\begin{proof}[Proof of Theorem~\ref{bs}] Assume that the 
tuple $f=(f^1,\ldots, f^m)$ is given. Choose disjoint isomorphic
bundles $E_j\simeq \C^m$ with isomorphic bases $e_{ji}$, and let
$f_j=\sum_{i=1}^m f^i e_{ji}^*$. Outside $Z=\{f=0\}$ we have
 $\sigma_j=\sum_1^m \sigma^i e_{ji}$. Now
$\dbar \sigma^i$ are linearly dependent, since
$\sum_1^m f^i\dbar\sigma^i=\dbar\sum_1^m f^i\sigma^i=\dbar 1=0$.
Thus the form  $u_k$ must vanish if $k-1>m-1$, and therefore  $R_k$ vanishes unless $k\le\min(n,m)$. 
Since  $|f_j|=|f|$,  locally in the resolution, we have
$$
R_k=\dbar|f|^{2r\lambda}\w
\frac{\beta}{(f^0)^{k+r-1}}\Big|_{\lambda=0},
$$
and hence it is annihilated by $\phi$ if $|\phi|\le C |f|^{\min(m,n)+r-1}$.
\end{proof}

It remains to consider the case when the $f_j$ together define a complete intersection. 
The proof is very much inspired by similar proofs in \cite{W}.

\begin{proof}[Proof of Theorem~\ref{ma2}]
We now assume that
$\codim\{f_1=\cdots =f_r=0\}=m_1+\cdots +m_r$. In particular, $m_1+\cdots +m_r\le n$.
Let  $\xi$ be  a test form times $\phi$. If the support is small enough,
after a resolution of singularities and further localization, 
$R.\xi$ becomes  a sum of terms, the worst of which are like
$$
\int \dbar\big(|f^0_1|^{2\lambda}\cdots |f^0_r|^{2\lambda}\big)\w
\frac{s_1'\w\ldots\w s_r'\w(\dbar s_1')^{m_1-1}\w\ldots\w(\dbar s'_r)^{m_r-1}\w\tilde\xi\rho}
{(f_1^0)^{m_1}\cdots (f_r^0)^{m_r}}\Big|_{\lambda=0},
$$
where $\tilde\xi $ is the pull-back of $\xi$ and 
$\rho$ is a cut-off function in the resolution.  We may assume that
 each $f^0_j$ is a monomial times a non-vanishing factor
in a  local coordinate system $\tau_k$. 
Let $\tau$ be one of the coordinate factors in, say,  $f_1$ (with  order $\ell$), and consider the integral
that appears when $\dbar$ falls on $|\tau^\ell|^{2\lambda}$. 
If $\tau$ does not occur  in any other $f^0_j$, then the assumption $|\phi|\le C|f_1|^{m_1}$
implies that $\tilde\phi$ is divisible by $\tau^{\ell m_1}$. Hence $\tilde\phi$ and therefore 
also $\tilde\xi$ 
annihilates the singularity as before, so that the integral vanishes.
We now claim that if,  on the other hand,  $\tau$ occurs in some of the other factors,
then the integral vanishes because of the complete intersection assumption.
Thus let us assume that $\tau$ occurs in $f^0_2,\ldots, f^0_k$ but not in
$f^0_{k+1},\ldots, f^0_r$.
The forms  $s_j=|f_j|^2\sigma_j$ are smooth  and,  moreover,
$$
\tilde\gamma=\frac{s'_{k+1}\w\ldots\w s'_r\w(\dbar s'_{k+1})^{m_{k+1}-1}\w\ldots
\w(\dbar s'_r)^{m_r-1}\w\tilde\xi}{(f_{k+1}^0)^{m_{k+1}}\cdots (f_r^0)^{m_r}}
$$
is the pull-back of 
$$
\gamma=
\frac{s_{k+1}\w\ldots\w s_r\w(\dbar s_{k+1})^{m_{k+1}-1}\w\ldots
\w(\dbar s_r)^{m_r-1}\w\xi}{|f_{k+1}|^{2m_{k+1}}\cdots |f_r^0|^{2m_r}}.
$$
Since the form $\gamma$ has codegree $1+(m_1-1)+\cdots +(m_k-1)$ in $d\bar z$, which is strictly
less than $m_1+\cdots +m_k=\codim\{f_1=\cdots =f_k=0\}$, the anti-holomorphic factor
of the denominator vanishes on $\{f_1=\cdots =f_k=0\}$. Therefore,
each term of its pull-back vanishes  where $\tau=0$,  so it must contain either
a factor $\bar\tau$ or $d\bar\tau$.  
However, because of the assumption, the (pull-back) of the 
denominator contains no factor $\bar\tau$,  so each term
of $\tilde\gamma$ will contain $\bar\tau$ or $d\bar\tau$.
Therefore, the integral that appears when $\dbar$ falls on $|\tau|^{2\lambda\ell}$
will vanish when $\lambda=0$.
\end{proof}

\section{Explicit integral representation}

We are now going to supply explicit proofs of Theorems~\ref{bs} to \ref{ma3}.
Since all of them are  local, we may assume, using the notation from
Section~\ref{curr},  that $f\colon E\to Q$ and the function $\phi$ are holomorphic
in a convex neighborhood $X$ of the closure of the unit
ball $\B$ in $\C^n$. Moreover, we fix global holomorphic frames $e_j$ and $\epsilon_k$
for $E$ and $Q$ respectively, and use the trivial metric with respect to these frames.

To give a hint of how the formulas are built up, first  suppose that $f$ is a function in the unit
disk with no zeros on the unit circle, i.e., $n=r=m=1$.
The construction of the representation  \eqref{deco} is  a generalization of the simple one-variable
formula
\begin{equation}\label{buss}
\phi(z)=f(z)\int_{|\zeta|=1}\frac{1}{f(\zeta)}\frac{d\zeta}{\zeta-z}\phi(\zeta)+
\frac{1}{2\pi i}\int_{|\zeta|<1}\dbar\frac{1}{f}\w h(\zeta,z) \phi(\zeta),
\end{equation}
where $h=(f(\zeta)-f(z))d\zeta/(\zeta-z)2\pi i$, 
which follows directly from Cauchy's integral formula. Notice that the second term
vanishes as soon as $\phi$ annihilates the residue $R=\dbar(1/f)$.
Moreover, for an arbitrary holomorphic function $\phi$, this term interpolates
$\phi$ at each  zero of $f$ up to the order of the zero. 
If the order is one this follows immediately from  the simple observation that 
$\dbar(1/f)\w df/2\pi i$ is the  point mass at the zero.

Turning our attention to the general case, one can verify   that if we in 
$R=\dbar|f|^{2\lambda}\w u|_{\lambda=0}$ contract each $\sigma=\sum\sigma_j\otimes\epsilon_j$
with $\sum df_j\otimes\epsilon_j^*$, and contract $\lsigma$ with $dF$, we get a $d$-closed
$(*,*)$-current of order zero, which in some sense  generalizes the Lelong current over $Z$;
see \cite{A3} for the case when $r=1$.  The recipe to obtain a division-interpolation formula
like  \eqref{deco} (and  \eqref{buss}) is to replace the  differentials  by  Hefer forms, and finally
multiply by a Cauchy type form.  The idea is developed in a quite general setting in
\cite{A1} so we only sketch our  special situation here.

\smallskip
For fixed  $z\in X$, 
let    $\delta_{\eta}$ denote interior multiplication with the vector field 
$2\pi i \sum (\zeta_j-z_j)(\partial/\partial \zeta_j),$
and let $\nabla_{\eta}=\delta_{\eta}-\dbar$. 
Moreover, let  $\chi$ be a cutoff function  in $X$ that is identically
$1$ in a neighborhood  of $\overline{\B}$ and  let
$$
s(\zeta,z)=\frac{1}{2\pi i}\frac{\partial|\zeta|^2}{|\zeta|^2-\bar\zeta\cdot z}.
$$
Then for each $z\in \B$,  see  \cite{A1}, 
\begin{equation}\label{gdef}
g=\chi-\dbar\chi\w\frac{s}{\nabla_\eta s}=
\chi-\dbar\chi\w\sum_{k=1}^n \frac{1}{(2\pi i)^k}
\frac{\partial|\zeta|^2\w(\dbar\partial|\zeta|^2)^{k-1}}{(|\zeta|^2-\bar\zeta\cdot z)^k}
\end{equation}
a compactly supported $\nabla_\eta$-closed form such that  (lower indices denote bidegree)
$g_{0,0}(z)=1$. Moreover,  $g$ depends holomorphically on $z$.  

\smallskip
We then  choose holomorphic $(1,0)$-forms $h_j$ in $X$  (Hefer forms)  such that
$\delta_\eta h_j=f_j(\zeta)-f_j(z),$
and let $h=\sum_1^m  h_j\otimes \epsilon_j^*$.
We may also assume that $h_j$, and hence $h$, depend holomorphically on the parameter $z$.
Now
$\delta_h\colon E_{k+1}\to E_k$, for $k\ge 1$, and hence
$(\delta_h)_{k}\colon E_{k+1}\to E_1, \quad k\ge 0,$
if  $(\delta_h)_\ell=\delta_h^\ell/\ell!$.
It is easily seen that
\begin{equation}\label{hh}
\delta_\eta (\delta_h)_k=(\delta_h)_{k-1}\delta_f-\delta_{f(z)}(\delta_h)_{k-1}.
\end{equation}
So far $\delta_F$ has only acted on $(0,0)$-forms with values in $\Lambda^rE$.  We now
extend it to general $(p,q)$-forms, with the convention that one insert a minus sign
when $p+q$ is odd. Thus we let 
$$
\delta_F\alpha=(-1)^{(r+1)(\deg\alpha+1)}\delta_{f_r}\cdots\delta_{f_1}\otimes\delta_{\epsilon_1}
\cdots\delta_{\epsilon_r},
$$
where  $\deg\alpha$ is the degree of $\alpha$ in $\Lambda(E\otimes T^*(X))$.
With this convention $\delta_F$, as well as  $\delta_f$,
will anti-commute with $\dbar$ and $\delta_\eta$.
It is   possible, see \cite{A1} to find explicit holomorphic $(1,0)$-form-valued 
mappings $H^0_k\colon E_k\to \C$ (depending holomorphically on the parameter $z$), such that
\begin{equation}\label{HH}
\delta_\eta H^0_1=\delta_{F(\zeta)}-\delta_{F(z)},\quad  
\delta_\eta H^0_k= H^0_{k-1}\delta_{f(\zeta)}-\delta_{F(z)}(\delta_h)_{k-1},\   k\ge 2.
\end{equation}
If we define 
$$
H^1U=\sum_{k=1}^{\min(n+1,m-r+1)} (\delta_h)_{k-1} U_k,
$$
and
$$
H^0R=\sum_{k=1}^{\min(n,m-r+1)} H^0_k R_k.
$$
it follows from \eqref{hh}   and \eqref{HH} that  $g'=(\delta_{F(z)}H^1U+H^0R)\w g$ is $\nabla_\eta$-closed,
and $g'_{0,0}(z)=1$,  and therefore we have, see \cite{A1},

\begin{thm} 
If $\phi$ is holomorphic in $X$ and $g$ is the Cauchy form 
\eqref{gdef}, 
then we have the holomorphic decomposition 
\begin{equation}\label{repp}
\phi(z)=\delta_{F(z)}\int H^1U \w g\phi +\int H^0R\w g\phi, \quad z\in\B.
\end{equation}
\end{thm}

In particular, 
$
\Psi(z)=\int H^1U\w g\phi
$
is an explicit solution to  $\delta_{F(z)}\Psi=\phi$  if $R\phi=0$.
We now consider this solution in more detail.
In view of \eqref{snorre} and  \eqref{udef0} we have, outside $Z$, that
\begin{multline*}
(\delta_h)_{k-1}u_k= \\
\sum_{|\alpha|=k-1}
(\delta_{h_1})_{\alpha_1}\cdots(\delta_{h_r})_{\alpha_r}\big[
\sigma_1\w\ldots\w\sigma_r\w(\dbar\sigma_1)^{\alpha_1}\w\ldots\w(\dbar\sigma_r)^{\alpha_r}\big]
\otimes\epsilon^*.
\end{multline*}
Moreover, since we have the trivial metric, 
$$
\sigma_j=\sum_{i=1}^m\sigma_{ij}e_j,  \quad j=1,\ldots,r,
$$
are just the columns in the matrix $f^*(ff^*)^{-1}$.
Suppressing the non-vanishing section $\epsilon^*$, we have

\begin{cor}
Let $f$ be a generically surjective holomorphic $r\times m$-matrix in $X$
with rows $f_j$, 
considered as sections of the trivial bundle $E^*$,
and assume that the hypothesis of  Theorem~\ref{ma3} is fulfilled.
Then
$$
\psi(z)=\int H^1U \w g\phi
$$ 
is an explicit solution to $\delta_{F(z)}\psi(z)=\delta_{f_1(z)}\cdots\delta_{f_r(z)}\psi(z)=\phi(z)$
in $\B$, 
where $H^1U\phi$ is the value at $\lambda=0$ of (the analytic continuation of)
\begin{multline}\label{polly}
|f|^{2\lambda}\sum_{k=1}^{\min(n+1,m-r+1)}\\
\sum_{|\alpha|=k-1} 
(\delta_{h_1})_{\alpha_1}\cdots(\delta_{h_r})_{\alpha_r}\big[
\sigma_1\w\ldots\w\sigma_r\w(\dbar\sigma_1)^{\alpha_1}\w\ldots\w(\dbar\sigma_r)^{\alpha_r}\big] 
\phi . 
\end{multline} 
If  $m-r+1\le n$, then $H^1U\phi$ is locally integrable, and the value at $\lambda=0$ exists
in the ordinary sense.
\end{cor}

\begin{proof}It remains to verify the claim about the local integrability.
In fact, after a resolution of singularities, cf., \eqref{using}, it follows that
$U_k\phi$ is locally integrable if  $|\phi|\lesssim |F|^{k}$. 
If $m-r+1\le n$, then the sum terminates at $k=m-r+1$, and therefore the current
is locally integrable; otherwise the worst term is like  $U_{n+1}\phi$, 
and it will not be locally integrable in general.
\end{proof}

If all the $f_j$ take values in different bundles $E_j^*$
and $E=\oplus E_j$,  then
we can simplify the expression for $H^1U$ further. In this case, cf., Section~\ref{poo}, 
$$
\sigma_j=\sum_{i=1}^{m_j} \frac{\bar f_j^i}{|f_j|^2}e_{ij},  \quad j=1,\ldots,r.
$$
Moreover, with natural choices of  Hefer forms $h_j$,  $\delta_{h_j}$ will vanish
on forms with values in $E_k$ for $k\neq j$, and hence we get 

\begin{cor}
Let  $f_j$ be  $m_j$-tuples of functions, considered as sections of the trivial bundles
$E_j^*$ over $X$.  If the conditions of Theorem~\ref{ma1} or \ref{ma2} are fulfilled,
or if all $f_j$ are equal to some fixed $m$-tuple $f$, and the condition in Theorem~\ref{bs}
is fulfilled, then 
$$
\psi(z)=\int H^1U\phi \w g
$$ 
is an explicit solution to $\delta_{f_1(z)}\cdots\delta_{f_r(z)}\psi(z)=\phi(z)$
in $\B$, 
where $H^1U\phi$ is the value at $\lambda=0$ of (the analytic continuation of)
\begin{multline}\label{polly2}
|f|^{2\lambda}\sum_{k=1}^{n+1}
\sum_{|\alpha|=k-1} 
(\delta_{h_1})_{\alpha_1}[\sigma_1\w(\dbar\sigma_1)^{\alpha_1}]\w
\ldots\w
(\delta_{h_r})_{\alpha_r}[\sigma_r\w(\dbar\sigma_r)^{\alpha_r}]
\phi,
\end{multline} 
and  $N=\min(n+1,m-r+1)$. 
\end{cor}

In the case of Theorems~\ref{ma1} and \ref{ma2}, only terms such that 
$\alpha_j\le m_j$  actually occur.
In the case of  Theorem~\ref{bs} we have only terms such that $k\le m$.

\bigskip

The division formulas discussed here are different from Berndtsson's classical
formulas, \cite{BB}. As already mentioned, an explicit  proof of
Theorem~\ref{bs} in the case $r=1$,  based on Berndtsson's division formula,
already appeared in \cite{BGVY} (see the proof of Theorem~3.25), and 
the  case with general $r$ follows in essentially the same way, see, \cite{EG}. 
However, we see no way to prove any of the variations  discussed
in this paper by classical Berndtsson type formulas.

\def\listing#1#2#3{{\sc #1}:\ {\it #2},\ #3.}

\end{document}